\documentclass[10 pt]{article}
\usepackage[pagebackref]{hyperref}
\hypersetup{
colorlinks = true,
linkcolor = blue,
citecolor = blue}
\usepackage[utf8]{inputenc}
\usepackage[T1]{fontenc}

\usepackage[]{amsfonts}
\usepackage{amssymb}
\usepackage{amsmath,amsthm}
\def\couleur(#1 #2 #3)
	{
\special{ps::[end]}}

\def\bx#1{\setbox1=\hbox{\kern3pt{#1}\kern3pt}			
 \dimen1=\ht1 \advance\dimen1 by 3pt \dimen2=\dp1 \advance\dimen2 by 3pt
 \setbox1=\hbox{\vrule height\dimen1 depth\dimen2\box1\vrule}%
 \setbox1=\vbox{\hrule\box1\hrule}%
 \advance\dimen1 by .4pt \ht1=\dimen1
 \advance\dimen2 by .4pt \dp1=\dimen2 \box1\relax}

\def\wbb#1{\kern#1em}
\def\vci{\vrule  width.02em height1.47ex depth-.0ex}		
\def\11{{\rm\wbb{.2}\vci\wbb{-.37}1}}

\def\underset#1#2{\mathrel{\mathop{\kern0pt #2}\limits_{#1}}}

\def\overset#1#2{\mathrel{\mathop{\kern0pt #2}\limits^{#1}}}

\def\Supp{\mathop{\rm Supp}\nolimits}

\newtheorem{thm}{Theorem}[section]
\newtheorem{lem}[thm]{Lemma}
\newtheorem{prop}[thm]{Proposition}
\newtheorem{cor}[thm]{Corollary}
\newtheorem{defin}[thm]{Definition}
\newtheorem{rem}[thm]{Remark}

\begin{document}

\title{An Andreotti-Grauert theorem with $L^{r}$ estimates.}

\author{Eric Amar}

\date{}
\maketitle
 \renewcommand{\abstractname}{Abstract}

\begin{abstract}
By a theorem of Andreotti and Grauert if $\omega $ is a $(p,q)$
 current, $q<n,$ in a Stein manifold, $\bar \partial $ closed
 and with compact support, then there is a solution $u$ to $\bar
 \partial u=\omega $ still with compact support. The main result
 of this work is to show that if moreover $\omega \in L^{r}(dm),$
 where $m$ is a suitable "Lebesgue" measure on the Stein manifold,
 then we have a solution $u$ with compact support \emph{and}
  in $L^{r}(dm).$ We prove it by estimates in $L^{r}$ spaces with weights.\par 

\end{abstract}
\ \par 

\tableofcontents
\ \par 
\ \par 

\section{Introduction.}
\quad \quad 	Let $\omega $ be a $\bar \partial $ closed $(p,q)$ form in ${\mathbb{C}}^{n}$
 with compact support $K:=\Supp \omega $ and such that $\omega
 \in L^{r}({\mathbb{C}}^{n}),$ the Lebesgue space in ${\mathbb{C}}^{n}.$
 Setting $K$ in a ball ${\mathbb{B}}:=B(0,\ R)$ with $R$ big
 enough, we know, by a theorem of Ovrelid~\cite{Ovrelid71}, that
 we have a $(p,q-1)$ form $u\in L^{r}({\mathbb{B}})$ such that
 $\bar \partial u=\omega .$ On the other hand we also know, at
 least when $q<n,$ that there is a current $v$ with compact support
 such that $\bar \partial v=\omega ,$ by a theorem of Andreotti-Grauert~\cite{AndreGrauert62}.\
 \par 
So a natural question is: may we have a solution $u$ of $\bar
 \partial u=\omega $ with compact support \emph{and} in $L^{r}({\mathbb{C}}^{n})\
 ?$\ \par 
\ \par 
\quad \quad 	There is a work by H. Skoda~\cite{zeroSkoda} who proved such
 a result. Let $\displaystyle \Omega $ be a strictly pseudo-convex
 bounded domain  in ${\mathbb{C}}^{n}$ with smooth boundary then
 in~\cite[Corollaire p. 295]{zeroSkoda}, H. Skoda proved that
 if $\displaystyle f$ is a $\displaystyle (p,q)$-form with measure
 coefficients, $q<n,\ \bar \partial $ closed and with compact
 support in $\Omega ,$ then there is a solution $U$ to the equation
 $\bar \partial U=f$ such that $\displaystyle \ {\left\Vert{U}\right\Vert}_{L^{r}(\Omega
 )}\leq C(\Omega ,r){\left\Vert{f}\right\Vert}_{1},$ for any
 $r$ such that$\ 1<r<\frac{2n+2}{2n-1}$ and $U$ has zero boundary
 values in the sense of Stokes formula. This means that essentially
 $U$ has compact support and, because $\displaystyle \Omega $
 is bounded, $\ {\left\Vert{f}\right\Vert}_{1}\lesssim {\left\Vert{f}\right\Vert}_{L^{r}(\Omega
 )}.$ So he got the answer for $\displaystyle \Omega $ strictly
 pseudo-convex and $1<r<\frac{2n+2}{2n-1}.$\ \par 
\ \par 
\quad \quad 	We answered this question by the affirmative for any $\displaystyle
 r\in \lbrack 1,\infty \rbrack $ in a joint work with S. Mongodi~\cite{AmMong12}
 linearly by the "method of coronas". This method asks for extra
 $L^{r}$ conditions on derivatives of coefficients of $\omega
 ,$ when $q<n$; we shall denote the set of $\omega $ verifying
 these conditions ${\mathcal{W}}^{r}_{q}(\Omega ),$ as in~\cite{AmMong12}.\
 \par 
\ \par 
\quad \quad 	The aim of this work is to extend this result to Stein manifolds
 and get rid of the extra $L^{r}$ conditions ${\mathcal{W}}^{r}_{q}(\Omega
 ).$ For it we use a completely different approach inspired by
 the Serre duality~\cite{Serre55}. Because Hahn Banach theorem
 is used, this method is no longer constructive as in~\cite{AmMong12}.\ \par 
\ \par 
\quad The basic notion we shall use here is the following.\ \par 

\begin{defin}
~\label{AG16}Let $X$ be a complex manifold equipped with a Borel
 $\sigma $-finite measure $dm$ and $\Omega $ a domain in $X$;
 let $\displaystyle r\in \lbrack 1,\ \infty \rbrack ,$ we shall
 say that $\Omega $ is $r$ {\bf regular} if for any $p,q\in \lbrace
 0,...,n\rbrace ,\ q\geq 1,$ there is a constant $C=C_{p,q}(\Omega
 )$ such that for any $(p,q)$ form $\omega ,\ \bar \partial $
 closed in $\Omega $ and in $L^{r}(\Omega ,dm)$ there is a $(p,q-1)$
 form $u\in L^{r}(\Omega ,dm)$ such that $\bar \partial u=\omega
 $ and $\ {\left\Vert{u}\right\Vert}_{L^{r}(\Omega )}\leq C{\left\Vert{\omega
 }\right\Vert}_{L^{r}(\Omega )}.$\par 
\quad \quad 	We shall say that $\Omega $ is {\bf weakly }$r$ {\bf regular}
 if for any compact set $K\Subset \Omega $ there are $3$ open
 sets $\Omega _{1},\Omega _{2},\Omega _{3}$ such that $K\Subset
 \Omega _{3}\subset \Omega _{2}\subset \Omega _{1}\subset \Omega
 _{0}:=\Omega $ and $3$ constants $C_{1},C_{2},C_{3}$ such that:\par 
\vspace{5pt} \ \ \ \ \ \ \ \ \ \ \ \  $\displaystyle \forall
 j=0,1,2,\ \forall p,q\in \lbrace 0,...,n\rbrace ,\ q\geq 1,\
 \forall \omega \in L_{p,q}^{r}(\Omega _{j},dm),\ \bar \partial
 \omega =0,$\vspace{5pt}  \par 
\vspace{5pt} \ \ \ \ \ \ \ \ \ \ \ \              $\displaystyle
 \exists u\in L_{p,q-1}^{r}(\Omega _{j+1},dm),\ \bar \partial
 u=\omega $\vspace{5pt}  \par 
and $\ {\left\Vert{u}\right\Vert}_{L^{r}(\Omega _{j+1})}\leq
 C_{j+1}{\left\Vert{\omega }\right\Vert}_{L^{r}(\Omega _{j})}.$\par 
\quad \quad 	I.e. we have a $3$ steps chain of resolution.
\end{defin}
Of course the $r$ regularity implies the weak $r$ regularity,
 just taking $\Omega _{1}=\Omega _{2}=\Omega _{3}=\Omega .$\ \par 
\ \par 
\quad \quad 	Examples of $2$ regular domains are the bounded pseudo-convex
 domains by H\"ormander~\cite{Hormander73}.\ \par 
\quad \quad 	Examples of $r$ regular domains in ${\mathbb{C}}^{n}$ are the
 bounded strictly pseudo-convex (s.p.c.) domains with smooth
 boundary by Ovrelid~\cite{Ovrelid71}; the polydiscs in ${\mathbb{C}}^{n}$
 by Charpentier~\cite{dbarCharp}, finite transverse intersections
 of strictly pseudo-convex bounded domains in ${\mathbb{C}}^{n}$
 by Menini~\cite{Menini97}. A generalisation of the results by
 Menini was done in the nice work of Ma and Vassiliadou~\cite{MaVassiliadou00}:
 they treated also the case of intersection of $q$-convex sets.\ \par 
\ \par 
\quad Examples of $r$ regular domains in a Stein manifold are the strictly
 pseudo-convex domains with smooth boundary~\cite{AmarSt13}.
 (See the previous work for $(0,1)$ forms by N. Kerzman~\cite{Kerzman71}
 and for all $(p,q)$ forms by J-P. Demailly and C. Laurent~\cite[Remarque
 4, page 596]{DemaillyLaurent87}, but here the manifold has to
 be equipped with a metric with null curvature. See also~\cite{dBarIntersection17}
 for the case of intersection of $q$-convex sets in a Stein manifold).\ \par 
\ \par 
\quad Let $X$ be a Stein manifold and $\Omega $ a domain in $X,$ i.e.
 an open connected set in $X.$ 	Let ${\mathcal{H}}_{p}(\Omega
 )$ be the set of all $(p,\ 0)\ \bar \partial $ closed forms
 in $\Omega .$ If $p=0,\ {\mathcal{H}}_{0}(\Omega )={\mathcal{H}}(\Omega
 )$  is the set of holomorphic functions in $\Omega .$ If $p>0,$
 we have, in a chart $(\varphi ,U),\ h\in {\mathcal{H}}_{p}(\Omega
 )\Rightarrow h(z)=\sum_{\left\vert{J}\right\vert =p}{a_{J}(z)dz^{J}},$
 where $dz^{J}:=dz_{j_{1}}\wedge \cdot \cdot \cdot \wedge dz_{j_{p}}$
 and the functions $a_{J}(z)$ are holomorphic in $\varphi (U)\subset
 {\mathbb{C}}^{n}.$\ \par 
\quad We shall denote $L_{p,q}^{r,c}(\Omega )$ the set of $(p,q)$-forms
 in $L^{r}(\Omega )$ with compact support in $\Omega .$ \ \par 
We also use the notation $\displaystyle r'$ for the conjugate
 exponent of $r,$ i.e. $\ \frac{1}{r}+\frac{1}{r'}=1.$\ \par 
\ \par 
\quad Our main theorem is:\ \par 

\begin{thm}
Let $\Omega $ be a weakly $r'$ regular domain in a Stein manifold
 $X.$ Then there is a $C>0$ such that for any  $(p,q)$ form $\omega
 $ in $L^{r,c}(\Omega ),\ r>1$ with:\par 
\quad $\bullet $ if $1\leq q<n,\ \bar \partial \omega =0$;\par 
\quad $\bullet $ if $q=n,\ \forall V\subset \Omega ,\ \Supp \omega
 \subset V,\ \omega \perp {\mathcal{H}}_{n-p}(V)$; \par 
there is a $(p,q-1)$ form $u$ in $L^{r,c}(\Omega )$ such that
 $\bar \partial u=\omega $ as distributions and 	${\left\Vert{u}\right\Vert}_{L^{r}(\Omega
 )}\leq C{\left\Vert{\omega }\right\Vert}_{L^{r}(\Omega )}.$
\end{thm}
\quad The notion of $r$ regularity gives a good control of the support:
 if the support of the data $\omega $ is contained in $\Omega
 \backslash C$ where $\Omega $ is a weakly $r'$ regular domain
 and $C$ is a weakly $r$ regular domain, then the support of
 the solution $u$ is contained in $\Omega \backslash C',$ where
 $C'$ is any relatively compact domain in $C,$ provided that
 $q\geq 2.$ One may observe that $\Omega \backslash C$ is \emph{not}
 Stein in general even if $\Omega $ is.\ \par 
\quad There is also a result of this kind for $q=1,$ see section~\ref{AndreGrau1913}.\
 \par 
\quad In particular the support of the solution $u$ is contained in
 the intersection of all the weakly $r'$ regular domains containing
 the support of $\omega .$\ \par 
\ \par 
\quad \quad 	The idea is to solve $\bar \partial u=\omega $ in a space $L^{r}(\Omega
 )$ with a "big weight $\eta $ outside" of the support of $\omega
 $; this way we shall have a "small solution $u$ outside" of
 the support of $\omega .$ Then, using a sequence of such weights
 going to infinity "outside" of the support of $\omega ,$ we
 shall have a $u$ zero "outside of the support" of $\omega .$\ \par 
\ \par 
\quad Comparing to my previous work~\cite{AnGrauLr1} the results here
 are improved and the proofs are much simpler by a systematic
 use of the Hodge $*$ operator.\ \par 
\ \par 
\quad \quad 	I am indebted to G. Tomassini who started my interest in this
 subject on precisely this kind of questions and also to S. Mongodi
 for a lot of discussions during the preparation of our joint
 paper~\cite{AmMong12}.\ \par 
\quad \quad 	I want to thank C. Laurent for many instructive discussions
 on this subject.\ \par 
\quad Finally I also thank the referee for his/her careful reading
 of the manuscript and the incisive questions he/she asked.\ \par 

\section{Duality.}
\quad \quad 	We shall study a duality between currents inspired by the Serre
 duality~\cite{Serre55}.\ \par 
\quad Let $X$ be a complex manifold of dimension $n.$ We proceed now
 exactly as in H\"ormander ~\cite[p. 119]{Hormander73}, by introducing
 a hermitian metric on differential  forms locally equivalent
 to the usual one on any analytic coordinates system.\ \par 
We define the "Lebesgue measure" still as in H\"ormander's book~\cite[Section
 5.2]{Hormander73}: associated to this metric there is a volume
 measure $dm$ and we take it for the Lebesgue measure on $X.$
 Moreover, because $X$ is a complex manifold, it is canonically oriented.\ \par 

\subsection{Weighted $L^{r}$ spaces.}
\quad Let $\Omega $ be a domain in $X.$ We denote also $dm$ the volume
 form on $X.$\ \par 
We shall take the following notation from the book by C. Voisin~\cite{Voisin02}.\
 \par 
\quad To a $(p,q)$-form $\alpha $ on $\Omega $ we associate its Hodge
 $*\ (n-p,n-q)$-form $*\alpha .$ This gives us a pointwise scalar
 product and a pointwise modulus:\ \par 
\quad \quad \quad \begin{equation}  (\alpha ,\beta )dm:=\alpha \wedge {\overline{*\beta
 }};\ \ \left\vert{\alpha }\right\vert ^{2}dm:=\alpha \wedge
 {\overline{*\alpha ,}}\label{AG23}\end{equation}\ \par 
because $\displaystyle \alpha \wedge {\overline{*\beta }}$ is
 a $(n,n)$-form hence is a function time the volume form $dm.$\ \par 
\quad We are also given a scalar product ${\left\langle{\alpha ,\beta
 }\right\rangle}$ on $(p,q)$-forms such that$\displaystyle \int_{\Omega
 }{\left\vert{\alpha }\right\vert ^{2}dm}<\infty $ and the link
 between these notions is given by~\cite[Lemme 5.8, p. 119]{Voisin02}:\ \par 
\quad \quad \quad \begin{equation} {\left\langle{ \alpha ,\beta }\right\rangle}=\int_{\Omega
 }{\alpha \wedge {\overline{*\beta }}}.\label{AG18}\end{equation}\ \par 
We shall define now $L_{p,q}^{r}(\Omega )$ to be the set of $(p,q)-$forms
 $\alpha $ defined on $\Omega $ such that\ \par 
\vspace{5pt} \ \ \ \ \ \ \ \ \ \ \ \  $\displaystyle {\left\Vert{\alpha
 }\right\Vert}_{L_{p,q}^{r}(\Omega )}^{r}:=\int_{\Omega }{\left\vert{\alpha
 (z)}\right\vert ^{r}dm(z)}<\infty ,$\vspace{5pt}  \ \par 
where $\left\vert{\alpha }\right\vert $ is defined by~(\ref{AG23}).\ \par 

\begin{lem}
~\label{AG22}Let $\eta >0$ be a weight. If $u$ is a $(p,q)$-current
 defined on $(n-p,n-q)$-forms $\alpha $ in $L^{r'}(\Omega ,\eta
 )$ and such that\par 
\vspace{5pt} \ \ \ \ \ \ \ \ \ \ \ \  $\displaystyle \forall
 \alpha \in L^{r'}_{(n-p,n-q)}(\Omega ,\eta ),\ \left\vert{{\left\langle{u,\
 *\alpha }\right\rangle}}\right\vert \leq C{\left\Vert{\alpha
 }\right\Vert}_{L^{r'}(\Omega ,\eta )},$\vspace{5pt}  \par 
then $\ {\left\Vert{u}\right\Vert}_{L_{p,q}^{r}(\Omega ,\eta ^{1-r})}\leq C.$
\end{lem}
\quad \quad 	Proof.\ \par 
We use the classical trick: set $\tilde \alpha :=\eta ^{1/r'}\alpha
 ;\ \tilde u:=\frac{1}{\eta ^{1/r'}}u$ then we have\ \par 
\vspace{5pt} \ \ \ \ \ \ \ \ \ \ \ \  $\displaystyle {\left\langle{u,*\alpha
 }\right\rangle}=\int_{\Omega }{u\wedge {\overline{\alpha }}}=\int_{\Omega
 }{\tilde u\wedge {\overline{\tilde \alpha }}}={\left\langle{\tilde
 u,*\tilde \alpha }\right\rangle}$\vspace{5pt}  \ \par 
and ${\left\Vert{\tilde \alpha }\right\Vert}_{L^{r'}(\Omega )}={\left\Vert{\alpha
 }\right\Vert}_{L^{r'}(\Omega ,\eta )}.$\ \par 
\quad We notice that ${\left\Vert{\tilde \alpha }\right\Vert}_{L^{r'}(\Omega
 )}={\left\Vert{*\tilde \alpha }\right\Vert}_{L^{r'}(\Omega )}$
 because we have $(*\tilde \alpha ,*\tilde \alpha )dm=*\tilde
 \alpha \wedge {\overline{**\tilde \alpha }}$ but $**\tilde \alpha
 =(-1)^{(p+q)(2n-p-q)}\tilde \alpha ,$ by ~\cite[Lemma 5.5]{Voisin02},
 hence, because $\displaystyle (*\tilde \alpha ,*\tilde \alpha
 )$ is positive,  $(*\tilde \alpha ,*\tilde \alpha )=\left\vert{\tilde
 \alpha }\right\vert ^{2}.$\ \par 
By use of the duality $L^{r}_{p,q}(\Omega )-L^{r'}_{n-p,n-q}(\Omega
 ),$ done in Lemma~\ref{AG25}, we get\ \par 
\vspace{5pt} \ \ \ \ \ \ \ \ \ \ \ \  $\displaystyle {\left\Vert{\tilde
 u}\right\Vert}_{L_{p,q}^{r}(\Omega )}=\sup _{\alpha \in L_{n-p,n-q}^{r'}(\Omega
 ),\ \alpha \neq 0}\frac{\left\vert{{\left\langle{\tilde u,*\tilde
 \alpha }\right\rangle}}\right\vert }{{\left\Vert{\tilde \alpha
 }\right\Vert}_{L^{r'}(\Omega )}}.$\vspace{5pt}  \ \par 
But\ \par 
\vspace{5pt} \ \ \ \ \ \ \ \ \ \ \ \  $\displaystyle {\left\Vert{\tilde
 u}\right\Vert}_{L_{p,q}^{r}(\Omega )}^{r}:=\int_{\Omega }{\left\vert{u}\right\vert
 ^{r}\eta ^{-\frac{r}{r'}}dm}=\int_{\Omega }{\left\vert{u}\right\vert
 ^{r}\eta ^{1-r}dm}={\left\Vert{u}\right\Vert}_{L^{r}(\Omega
 ,\eta ^{1-r})}^{r}.$\vspace{5pt}  \ \par 
So we get\ \par 
\vspace{5pt} \ \ \ \ \ \ \ \ \ \ \ \  $\displaystyle {\left\Vert{u}\right\Vert}_{L_{p,q}^{r}(\Omega
 ,\eta ^{1-r})}=\sup _{*\alpha \in L_{p,q}^{r'}(\Omega ,\eta
 ),\ \alpha \neq 0}\frac{\left\vert{{\left\langle{u,*\alpha }\right\rangle}}\right\vert
 }{{\left\Vert{\alpha }\right\Vert}_{L^{r'}(\Omega ,\eta )}}.$\vspace{5pt}
  \ \par 
\quad The proof is complete. $\hfill\blacksquare $\ \par 
\ \par 
\quad It may seem strange that we have such an estimate when the dual
 of $L^{r'}(\Omega ,\eta )$ is  $L^{r}(\Omega ,\eta ),$ but the
 reason is, of course, that in the duality current-form there
 is no weights.\ \par 
\quad The point here is that when $\eta $ is small, $\eta ^{1-r}$ is
 big for $r>1.$\ \par 

\section{Solution of the $\bar \partial $ equation with compact support.}

\subsection{$r${\it  } regular domains.}
\quad As we have seen, examples of $r$ regular domains in Stein manifolds
 are the relatively compact s.p.c. domains with smooth boundary.\ \par 
\quad To prove that a Stein manifold $\Omega $ is weakly $r$ regular
 we shall need the following lemma.\ \par 

\begin{lem}
~\label{AAG26} Let $\Omega $ be a Stein manifold. Then it contains
 an exhaustive sequence of open relatively compact strictly pseudo-convexs
 sets $\lbrace D_{k}\rbrace _{k\in {\mathbb{N}}}$ with ${\mathcal{C}}^{\infty
 }$ smooth boundary.
\end{lem}
\quad Proof.\ \par 
For the case of $\Omega $ pseudo-convex in ${\mathbb{C}}^{n},$
 the proof was already done explicitely in the proof of~\cite[Theorem
 2.8.1, p. 86]{HenkinLeiterer84}.\ \par 
\quad By Theorem 5.1.6 of H\"ormander~\cite{Hormander73} there exists
 a ${\mathcal{C}}^{\infty }$ strictly plurisubharmonic (s.p.s.h.)
 exhausting function $\varphi $ for $\Omega .$ Take $K\Subset
 \Omega $ such that $d\varphi \neq 0$ on $K.$ Because $\varphi
 $ is s.p.s.h. then $K\neq \emptyset .$ Then we use the~\cite[Lemma
 2.12.2, p. 93]{HenkinLeiterer84}, to get:\ \par 
\quad \quad \quad $\forall \epsilon >0,\ \exists \rho _{\epsilon }\ s.p.s.h.\ {\mathcal{C}}^{\infty
 }$-function on $\Omega $ such that:\ \par 
\quad (i) $\varphi -\rho _{\epsilon }$ together with its first and
 second derivatives is less than $\epsilon $ on $\Omega .$\ \par 
\quad (ii) The set $\mathrm{C}\mathrm{r}\mathrm{i}\mathrm{t}(\rho _{\epsilon
 }):=\lbrace z\in \Omega ::d\rho _{\epsilon }(z)=0\rbrace $ is
 discrete in $\Omega .$ (In a formula, the notation $::$ means
 "such that".)\ \par 
\quad (iii) $\rho _{\epsilon }=\varphi $ on $K.$\ \par 
As stated in Lemma 2.12.2 if $\varphi \in {\mathcal{C}}^{2}$
 then $\rho _{\epsilon }\in {\mathcal{C}}^{2},$ but in fact the
 proof of this Lemma 2.12.2 gives $\rho _{\epsilon }=\varphi
 +\sum{\chi _{j}},$ where $\displaystyle \sum{\chi _{j}}$ is
 ${\mathcal{C}}^{\infty }.$ (see p. 93 in~\cite{HenkinLeiterer84}).
 Hence $\rho _{\epsilon }$ has the same ${\mathcal{C}}^{k}$ regularity
 than $\varphi .$\ \par 
\quad Fix $\epsilon >0,$ then the  function $\rho :=\rho _{\epsilon
 }$ is also a s.p.s.h. exhausting function for $\Omega ,$ because,
 from $-\epsilon \leq \varphi -\rho _{\epsilon }\leq \epsilon
 ,$ we get that, for any $\alpha \in {\mathbb{R}},$ \ \par 
\vspace{5pt} \ \ \ \ \ \ \ \ \ \ \ \  $\displaystyle \lbrace
 z\in \Omega ::\rho _{\epsilon }(z)<\alpha \rbrace \subset \lbrace
 z\in \Omega ::\varphi (z)<\epsilon +\alpha \rbrace $\vspace{5pt}  \ \par 
and, because $\varphi $ is an exhausting function, this set is
 relatively compact in $\Omega .$\ \par 
\ \par 
\quad Because the set of critical points of $\rho $ is discrete in
 $\Omega ,$ the same way as in the proof of~\cite[Theorem 2.8.1,
 p. 86]{HenkinLeiterer84}, we can find a sequence $c_{k}\in {\mathbb{R}},\
 c_{k}\rightarrow \infty ,$ such that	 	$\displaystyle D_{k}:=\lbrace
 z\in \Omega ::\rho (z)<c_{k}\rbrace $ make an exhaustive sequence
 of open relatively compact sets in $\Omega ,\ \partial \rho
 \neq 0$ on $\partial D_{k},$ hence $D_{k}$ is strictly pseudo-convex
 with ${\mathcal{C}}^{\infty }$ smooth boundary, and finally
 $D_{k}\nearrow \Omega .$ The proof is complete. $\hfill\blacksquare $\ \par 

\begin{prop}
A Stein manifold $\Omega $ is weakly $r$ regular.
\end{prop}
\quad Proof.\ \par 
By Lemma~\ref{AAG26} there is an exhaustive sequence of open
 relatively compact s.p.c. sets in $\Omega ,\ \lbrace D_{k}\rbrace
 _{k\in {\mathbb{N}}}$ with ${\mathcal{C}}^{\infty }$ smooth boundary.\ \par 
\quad Let $\omega \in L^{r}_{p,q}(\Omega ),\ \bar \partial \omega =0,$
 by~\cite{AmarSt13}, we can solve $\bar \partial u=\omega $ in
 $D_{k}$ with $u\in L^{r}_{p,q-1}(D_{k})$ and\ \par 
\vspace{5pt} \ \ \ \ \ \ \ \ \ \ \ \  $\displaystyle {\left\Vert{u}\right\Vert}_{L^{r}(D_{k})}\leq
 C_{k}{\left\Vert{\omega }\right\Vert}_{L^{r}(D_{k})}\leq C_{k}{\left\Vert{\omega
 }\right\Vert}_{L^{r}(\Omega )}.$\vspace{5pt}  \ \par 
Hence if $\Gamma $ is a compact set in $\Omega ,$ there is a
 $D_{k}$ such that $\Gamma \Subset D_{k}$ and we can take $\Omega
 _{1}=\Omega _{2}=\Omega _{3}=D_{k}.$\ \par 
\quad This proves the weak $r$ regularity of $\Omega .$ $\hfill\blacksquare $\ \par 

\subsection{The main result.}
\quad Let $X$ be a Stein manifold and $\Omega $ a domain in $X.$\ \par 
\quad In order to simplify notation, we set the pairing for $\alpha
 $ a $(p,q)$-form and $\beta $ a $(n-p,n-q)$-form: $\displaystyle
 \ll \alpha ,\beta \gg :=\int_{\Omega }{\alpha \wedge \beta }.$\ \par 
With this notation we also have ${\left\langle{\alpha ,\beta
 }\right\rangle}=\ll \alpha ,{\overline{*\beta }}\gg .$\ \par 
\ \par 
\quad \quad 	Let $\Omega $ be a weakly $r'$ regular domain in $X.$ We set
 $K:=\Supp \omega \Subset \Omega $ and, by the definition of
 the $r'$ weak regularity, we get $3$ open sets such that $K\Subset
 \Omega _{3}\subset \Omega _{2}\subset \Omega _{1}\subset \Omega
 _{0}=\Omega $ with: $\forall j=0,1,2,\ \forall p,q\in \lbrace
 0,...,n\rbrace ,\ q\geq 1,$ \ \par 
\vspace{5pt} \ \ \ \ \ \ \ \ \ \ \ \  $\displaystyle \forall
 \alpha \in L_{p,q}^{r}(\Omega _{j}),\ \bar \partial \alpha =0,\
 \exists \varphi \in L_{p,q-1}^{r}(\Omega _{j+1}),\ \bar \partial
 \varphi =\alpha .$\vspace{5pt}  \ \par 
\quad Set the weight $\eta =\eta _{\epsilon }:={\11}_{\Omega _{1}}(z)+\epsilon
 {\11}_{\Omega \backslash \Omega _{1}}(z)$ for a fixed $\epsilon >0.$\ \par 
\quad Let $\omega \in L_{p,q}^{r,c}(\Omega ).$ Suppose moreover that
 $\omega $ is such that $\bar \partial \omega =0$ if $1\leq q<n$
 and for any open $V\Subset \Omega ,\ \Supp \omega \Subset V$
 we have $\omega \perp {\mathcal{H}}_{n-p}(V)\iff \forall h\in
 {\mathcal{H}}_{n-p}(V),\ \ll \omega ,h\gg =0$ if $q=n.$\ \par 
\ \par 
We shall use the following lemma, with the previous notation:\ \par 

\begin{lem}
~\label{AAG27} Let ${\mathcal{E}}$ be the set of $(n-p,n-q+1)$
 forms $\alpha \in L^{r'}(\Omega ,\eta ),\ \bar \partial $ closed
 in $\Omega .$ Let us define $\displaystyle {\mathcal{L}}_{\omega
 }$ on ${\mathcal{E}}$ as follows:\par 
\vspace{5pt} \ \ \ \ \ \ \ \ \ \ \ \  $\displaystyle {\mathcal{L}}_{\omega
 }(\alpha ):=(-1)^{p+q-1}\ll \varphi ,\omega \gg ,$\vspace{5pt}  \par 
where $\varphi \in L^{r'}(\Omega _{1})$ is such that $\bar \partial
 \varphi =\alpha $ in $\Omega _{1}.$ Then the form ${\mathcal{L}}_{\omega
 }$ is well defined and linear.
\end{lem}
\quad \quad 	Proof.\ \par 
Because $\epsilon >0$ we have $\alpha \in L^{r'}(\Omega ,\eta
 )\Rightarrow \alpha \in L^{r'}(\Omega )$ and the weak $r'$ regularity
 of $\Omega $ gives a $\varphi \in L^{r'}(\Omega _{1})$ with
 $\bar \partial \varphi =\alpha $ in $\Omega _{1}.$\ \par 
Let us see that ${\mathcal{L}}_{\omega }$ is well defined.\ \par 
\ \par 
\quad $\bullet $	 Suppose first that $q<n.$\ \par 
In order for ${\mathcal{L}}_{\omega }$ to be well defined we need\ \par 
\vspace{5pt} \ \ \ \ \ \ \ \ \ \ \ \  $\displaystyle \forall
 \varphi ,\psi \in L^{r'}_{(n-p,n-q)}(\Omega _{1}),\ \bar \partial
 \varphi =\bar \partial \psi =\alpha \Rightarrow \ll \varphi
 ,\omega \gg =\ll \psi ,\omega \gg .$\vspace{5pt}  \ \par 
This is meaningful because $\omega \in L^{r,c}(\Omega ),\ r>1,\
 \Supp \omega \Subset \Omega _{1}.$\ \par 
Then we have $\bar \partial (\varphi -\psi )=0$ in $\Omega _{1},$
 hence, because $\Omega $ is weakly $r'$ regular, we can solve
 $\bar \partial $ in $L^{r'}(\Omega _{2})$:\ \par 
\vspace{5pt} \ \ \ \ \ \ \ \ \ \ \ \  $\displaystyle \exists
 \gamma \in L^{r'}_{(n-p,n-q-1)}(\Omega _{2})::\bar \partial
 \gamma =(\varphi -\psi ).$\vspace{5pt}  \ \par 
So $\ \ll \varphi -\psi ,\omega \gg =\ll \bar \partial \gamma
 ,\omega \gg =(-1)^{p+q-1}\ll \gamma ,\bar \partial \omega \gg
 =0$ because $\omega $ is compactly supported in $\Omega _{2}$
 and $\bar \partial $ closed.\ \par 
\quad Hence ${\mathcal{L}}_{\omega }$ is well defined in that case.\ \par 
\ \par 
\quad $\bullet $	 Suppose now that $q=n.$\ \par 
For $\varphi ,\ \psi \ (n-p,\ 0)$ forms in $\Omega _{1},$ such
 that $\bar \partial \varphi =\bar \partial \psi =\alpha ,$ we
 need to have $\ll \varphi ,\omega \gg =\ll \psi ,\omega \gg
 .$ But then $\bar \partial (\varphi -\psi )=0,$ which means
 that $h:=\varphi -\psi $ is a $\bar \partial $ closed $(n-p,\
 0)$ form, hence $h\in {\mathcal{H}}_{n-p}(\Omega _{1}).$ Taking
 $V=\Omega _{1}$ in the hypothesis $\omega \perp {\mathcal{H}}_{n-p}(V),$
 we get $\ll h,\omega \gg =0,$ and ${\mathcal{L}}_{\omega }$
 is also well defined in that case.\ \par 
\ \par 
\quad \quad 	It remains to see that ${\mathcal{L}}_{\omega }$ is linear.\ \par 
\quad $\bullet $ Suppose first that $q<n.$\ \par 
Let $\alpha =\alpha _{1}+\alpha _{2},$ with $\alpha _{j}\in L^{r'}(\Omega
 ,\eta ),\ \bar \partial \alpha _{j}=0,\ j=1,2$; we have $\alpha
 =\bar \partial \varphi ,\ \alpha _{1}=\bar \partial \varphi
 _{1}$ and $\alpha _{2}=\bar \partial \varphi _{2},$ with $\varphi
 ,\ \varphi _{1},\ \varphi _{2}$ in $\displaystyle L^{r'}(\Omega
 _{1})$ so, because $\bar \partial (\varphi -\varphi _{1}-\varphi
 _{2})=0,$ we have\ \par 
\vspace{5pt} \ \ \ \ \ \ \ \ \ \ \ \  $\displaystyle \ \varphi
 =\varphi _{1}+\varphi _{2}+\bar \partial \psi ,$\vspace{5pt}
   with $\psi $ in $L^{r'}(\Omega _{2}),$\vspace{5pt}  \ \par 
so\ \par 
\vspace{5pt} \ \ \ \ \ \ \ \ \ \ \ \  $\displaystyle {\mathcal{L}}_{\omega
 }(\alpha )=(-1)^{p+q-1}\ll \varphi ,\omega \gg =(-1)^{p+q-1}\ll
 \varphi _{1}+\varphi _{2}+\bar \partial \psi ,\omega \gg =$\vspace{5pt}
  \ \par 
\vspace{5pt} \ \ \ \ \ \ \ \ \ \ \ \ \ \ \ \ \ \ \ \ \ \ \ \
  $\displaystyle =\ {\mathcal{L}}_{\omega }(\alpha _{1})+{\mathcal{L}}_{\omega
 }(\alpha _{2})+(-1)^{p+q-1}\ll \bar \partial \psi ,\omega \gg
 ,$\vspace{5pt}  \ \par 
but again $\ll \bar \partial \psi ,\omega \gg =0,$ hence ${\mathcal{L}}_{\omega
 }(\alpha )={\mathcal{L}}_{\omega }(\alpha _{1})+{\mathcal{L}}_{\omega
 }(\alpha _{2}).$\ \par 
The same for $\alpha =\lambda \alpha _{1}.$ \ \par 
\ \par 
\quad $\bullet $	 Suppose now that $q=n.$\ \par 
We have\ \par 
\vspace{5pt} \ \ \ \ \ \ \ \ \ \ \ \  $\displaystyle \ {\mathcal{L}}_{\omega
 }(\alpha ):=(-1)^{p+n-1}\ll \varphi ,\omega \gg ,$\vspace{5pt}  \ \par 
where $\varphi \in L^{r'}(\Omega _{1})$ is such that $\bar \partial
 \varphi =\alpha $ in $\Omega _{1}.$\ \par 
Let $\alpha =\alpha _{1}+\alpha _{2},$ with $\alpha _{j}\in L^{r'}(\Omega
 ,\eta ),\ \bar \partial \alpha _{j}=0,\ j=1,2$; we have $\alpha
 =\bar \partial \varphi ,\ \alpha _{1}=\bar \partial \varphi
 _{1}$ and $\alpha _{2}=\bar \partial \varphi _{2},$ with $\varphi
 ,\ \varphi _{1},\ \varphi _{2}$ in $\displaystyle L^{r'}(\Omega
 _{1})$ so, because $\bar \partial (\varphi -\varphi _{1}-\varphi
 _{2})=0,$ we have $\displaystyle \varphi -\varphi _{1}-\varphi
 _{2}$ is a $(n-p,0)\ \bar \partial $-closed form, hence:\ \par 
\vspace{5pt} \ \ \ \ \ \ \ \ \ \ \ \  $\displaystyle \varphi
 =\varphi _{1}+\varphi _{2}+h,$ with $\displaystyle h\in {\mathcal{H}}_{n-p}(\Omega
 _{1}).$\vspace{5pt}  \ \par 
So\ \par 
\vspace{5pt} \ \ \ \ \ \ \ \ \ \ \ \  $\displaystyle {\mathcal{L}}_{\omega
 }(\alpha )=(-1)^{p+q-1}\ll \varphi ,\omega \gg =(-1)^{p+q-1}\ll
 \varphi _{1}+\varphi _{2}+h,\omega \gg =$\vspace{5pt}  \ \par 
\vspace{5pt} \ \ \ \ \ \ \ \ \ \ \ \ \ \ \ \ \ \ \ \ \ \ \ \
  $\displaystyle =\ {\mathcal{L}}_{\omega }(\alpha _{1})+{\mathcal{L}}_{\omega
 }(\alpha _{2})+(-1)^{p+q-1}\ll h,\omega \gg .$\vspace{5pt}  \ \par 
Taking $V=\Omega _{1}$ in the hypothesis $\omega \perp {\mathcal{H}}_{n-p}(V),$
 we get $\ll h,\omega \gg =0,$ hence ${\mathcal{L}}_{\omega }(\alpha
 )={\mathcal{L}}_{\omega }(\alpha _{1})+{\mathcal{L}}_{\omega
 }(\alpha _{2}).$\ \par 
The same for $\alpha =\lambda \alpha _{1}.$ The proof is complete.
 $\hfill\blacksquare $\ \par 

\begin{rem}
~\label{AndreGrau2815}If $\Omega $ is Stein, we can take the
 domain $\Omega _{1}$ to be s.p.c. with ${\mathcal{C}}^{\infty
 }$ smooth boundary, hence also Stein. So because $K:=\Supp \omega
 \subset \Omega _{1}\subset \Omega ,$ the $A(\Omega _{1})$ convex
 hull of $K,\ \hat K_{\Omega _{1}}$ is still in $\Omega _{1},$
 and any holomorphic function in $\Omega _{1}$ can be uniformly
 approximated on $\hat K_{\Omega _{1}}$ by holomorphic functions
 in $\Omega .$\par 
Then for $q=n$ instead of asking $\omega \perp {\mathcal{H}}_{n-p}(\Omega
 _{1})$ we need just $\omega \perp {\mathcal{H}}_{n-p}(\Omega ).$
\end{rem}

\begin{thm}
~\label{AndGrau84}Let $\Omega $ be a weakly $r'$ regular domain
 and  $\omega $ be a $(p,q)$ form in $L^{r,c}(\Omega ),\ r>1.$
 Suppose that $\omega $ is such that:\par 
\quad $\bullet $ if $1\leq q<n,\ \bar \partial \omega =0$; \par 
\quad $\bullet $ if $q=n,\ \forall V\subset \Omega ,\ \Supp \omega
 \subset V,\ \omega \perp {\mathcal{H}}_{n-p}(V)$\!\!\!\! .\par 
Then there is a $C>0$ and a $(p,q-1)$ form $u$ in $L^{r,c}(\Omega
 )$ such that $\bar \partial u=\omega $ as distributions and
 	$\displaystyle {\left\Vert{u}\right\Vert}_{L^{r}(\Omega )}\leq
 C{\left\Vert{\omega }\right\Vert}_{L^{r}(\Omega )}.$
\end{thm}
\quad \quad 	Proof.\ \par 
Because $\Omega $ is weakly $r'$ regular there is a $\Omega _{1}\subset
 \Omega ,\ \Omega _{1}\supset \Supp \omega $ such that\ \par 
\vspace{5pt} \ \ \ \ \ \ \ \ \ \ \ \  $\displaystyle \forall
 \alpha \in L^{r'}(\Omega ),\bar \partial \alpha =0,\ \exists
 \varphi \in L^{r'}(\Omega _{1})::\bar \partial \varphi =\alpha
 ,\ {\left\Vert{\varphi }\right\Vert}_{L^{r'}(\Omega _{1})}\leq
 C_{1}{\left\Vert{\alpha }\right\Vert}_{L^{r'}(\Omega )}.$\vspace{5pt}  \ \par 
There is a $\Omega _{2}$ such that $\Supp \omega \Subset \Omega
 _{2}\subset \Omega _{1}\subset \Omega $ with the same properties
 as $\Omega _{1}.$\ \par 
Let us consider the weight $\eta =\eta _{\epsilon }:={\11}_{\Omega
 _{1}}(z)+\epsilon {\11}_{\Omega \backslash \Omega _{1}}(z)$
 for a fixed $\epsilon >0$ and the form ${\mathcal{L}}_{\omega
 }$ defined in Lemma~\ref{AAG27}. By Lemma~\ref{AAG27} we have
 that ${\mathcal{L}}_{\omega }$ is a linear form on $(n-p,n-q+1)$-forms
 $\alpha \in L^{r'}(\Omega ,\eta ),\ \bar \partial $ closed in $\Omega .$\ \par 
\quad \quad 	If $\alpha $ is a  $(n-p,n-q+1)$-form in $L^{r'}(\Omega ,\eta
 ),$ then $\alpha $ is in $L^{r'}(\Omega )$ because $\epsilon >0.$\ \par 
The weak $r'$ regularity of $\Omega $ gives that there is a $\varphi
 \in L^{r'}(\Omega _{1})::\bar \partial \varphi =\alpha $ which
 can be used to define ${\mathcal{L}}_{\omega }(\alpha ).$\ \par 
\quad \quad 	We have also that $\alpha \in L^{r'}(\Omega _{1}),\ \bar \partial
 \alpha =0$ in $\Omega _{1},$ hence, still with the weak $r'$
 regularity of $\Omega ,$ we have\ \par 
\vspace{5pt} \ \ \ \ \ \ \ \ \ \ \ \  $\displaystyle \exists
 \psi \in L^{r'}(\Omega _{2})::\bar \partial \psi =\alpha ,\
 {\left\Vert{\psi }\right\Vert}_{L^{r'}(\Omega _{2})}\leq C_{2}{\left\Vert{\alpha
 }\right\Vert}_{L^{r'}(\Omega _{1})}.$\vspace{5pt}  \ \par 
\quad $\bullet $ For $q<n,$ we have $\bar \partial (\varphi -\psi )=\alpha
 -\alpha =0$ on $\Omega _{2}$ and, by the weak $r'$ regularity
 of $\Omega ,$ there is a $\Omega _{3}\subset \Omega _{2},\ $
 such that $\Supp \omega \subset \Omega _{3}\subset \Omega _{2},$
 and a $\gamma \in L^{r'}(\Omega _{3}),\ \bar \partial \gamma
 =\varphi -\psi $ in $\Omega _{3}.$ So we get\ \par 
\vspace{5pt} \ \ \ \ \ \ \ \ \ \ \ \  $\displaystyle \ll \varphi
 -\psi ,\omega \gg =\ll \bar \partial \gamma ,\omega \gg =(-1)^{p+q-1}\
 \ll \gamma ,\bar \partial \omega \gg =0,$\vspace{5pt}  \ \par 
this is meaningful  because $\Supp \omega \subset \Omega _{3}.$\ \par 
Hence	 	$\displaystyle {\mathcal{L}}_{\omega }(\alpha )=\ll \varphi
 ,\omega \gg =\ll \psi ,\omega \gg .$\ \par 
\ \par 
\quad $\bullet $	 For $q=n,$ we still have $\bar \partial (\varphi
 -\psi )=\alpha -\alpha =0$ on $\Omega _{2},$ hence $\displaystyle
 \varphi -\psi \in {\mathcal{H}}_{p}(\Omega _{2})$; this time
 we choose $V=\Omega _{2}$ and the assumption gives $\ \ll \varphi
 -\psi ,\omega \gg =0$ hence again ${\mathcal{L}}_{\omega }(\alpha
 )=\ll \varphi ,\omega \gg =\ll \psi ,\omega \gg .$\ \par 
\ \par 
\quad \quad 	In any cases, by H\"older inequalities done in Lemma~\ref{AG24},\ \par 
\vspace{5pt} \ \ \ \ \ \ \ \ \ \ \ \  $\displaystyle \left\vert{{\mathcal{L}}_{\omega
 }(\alpha )}\right\vert \leq {\left\Vert{\omega }\right\Vert}_{L^{r}(\Omega
 _{1})}{\left\Vert{\psi }\right\Vert}_{L^{r'}(\Omega _{2})}\leq
 {\left\Vert{\omega }\right\Vert}_{L^{r}(\Omega )}{\left\Vert{\psi
 }\right\Vert}_{L^{r'}(\Omega _{2})}.$\vspace{5pt}  \ \par 
But, by the weak $r'$ regularity of $\Omega ,$ there is a constant
 $C_{2}$ such that\ \par 
\vspace{5pt} \ \ \ \ \ \ \ \ \ \ \ \  $\displaystyle {\left\Vert{\psi
 }\right\Vert}_{L^{r'}(\Omega _{2})}\leq C_{2}{\left\Vert{\alpha
 }\right\Vert}_{L^{r'}(\Omega _{1})}.$\vspace{5pt}  \ \par 
Of course we have\ \par 
\vspace{5pt} \ \ \ \ \ \ \ \ \ \ \ \  $\displaystyle {\left\Vert{\alpha
 }\right\Vert}_{L^{r'}(\Omega _{1})}\leq {\left\Vert{\alpha }\right\Vert}_{L^{r'}(\Omega
 ,\ \eta )}$\vspace{5pt}  \ \par 
because $\eta =1$ on $\Omega _{1},$ hence\ \par 
\vspace{5pt} \ \ \ \ \ \ \ \ \ \ \ \  $\displaystyle \ \left\vert{{\mathcal{L}}_{\omega
 }(\alpha )}\right\vert \leq C_{2}{\left\Vert{\omega }\right\Vert}_{L^{r}(\Omega
 )}{\left\Vert{\alpha }\right\Vert}_{L^{r'}(\Omega ,\ \eta )}.$\vspace{5pt}
  \ \par 
\quad \quad 	So we have that the norm of ${\mathcal{L}}_{\omega }$ is bounded
 on the subspace of $\bar \partial $ closed forms in $L^{r'}(\Omega
 ,\eta )$ by $C{\left\Vert{\omega }\right\Vert}_{L^{r}(\Omega
 )}$ which is \emph{independent} of $\epsilon .$\ \par 
\quad \quad 	We apply the Hahn-Banach theorem to extend ${\mathcal{L}}_{\omega
 }$ with the \emph{same} norm to \emph{all}  $(n-p,n-q+1)$ forms
 in $L^{r'}(\Omega ,\ \eta ).$ As in the Serre Duality Theorem~\cite[p.
 20]{Serre55}, this is one of the major ingredients in the proof.\ \par 
\quad \quad 	This means, by the definition of currents, that there is a $(p,q-1)$
 current $u$ which represents the extended form ${\mathcal{L}}_{\omega
 }$: ${\mathcal{L}}_{\omega }(\alpha )=\ll \alpha ,u\gg .$ So
 if $\alpha :=\bar \partial \varphi $ with $\varphi \in {\mathcal{C}}^{\infty
 }_{c}(\Omega ),$ we get\ \par 
\vspace{5pt} \ \ \ \ \ \ \ \ \ \ \ \  $\displaystyle {\mathcal{L}}(\alpha
 )=\ll \alpha ,u\gg =\ll \bar \partial \varphi ,u\gg =(-1)^{p+q-1}\
 \ll \varphi ,\omega \gg $\vspace{5pt}  \ \par 
hence $\bar \partial u=\omega $  as distributions because $\varphi
 $ is compactly supported. And we have:\ \par 
\vspace{5pt} \ \ \ \ \ \ \ \ \ \ \ \  $\displaystyle \sup _{\alpha
 \in L^{r'}(\Omega ,\eta ),\ {\left\Vert{\alpha }\right\Vert}=1}\
 \left\vert{\ll \alpha ,u\gg }\right\vert \leq C{\left\Vert{\omega
 }\right\Vert}_{L^{r}(\Omega )}.$\vspace{5pt}  \ \par 
By lemma~\ref{AG22} with the weight $\eta ,$ this implies\ \par 
\vspace{5pt} \ \ \ \ \ \ \ \ \ \ \ \  $\displaystyle {\left\Vert{u}\right\Vert}_{L^{r}(\Omega
 ,\eta ^{1-r})}\leq C{\left\Vert{\omega }\right\Vert}_{L^{r}(\Omega
 )}$\vspace{5pt}  \ \par 
because $\left\vert{\ll \alpha ,u\gg }\right\vert =\left\vert{{\left\langle{\alpha
 ,{\overline{*u}}}\right\rangle}}\right\vert $ and, as already
 seen, ${\left\Vert{u}\right\Vert}_{L^{r}(\Omega ,\eta ^{1-r})}={\left\Vert{*u}\right\Vert}_{L^{r}(\Omega
 ,\eta ^{1-r})}={\left\Vert{{\overline{*u}}}\right\Vert}_{L^{r}(\Omega
 ,\eta ^{1-r})}.$\ \par 
In particular $\ {\left\Vert{u}\right\Vert}_{L^{r}(\Omega )}\leq
 C{\left\Vert{\omega }\right\Vert}_{L^{r}(\Omega )}$ because
 with $\epsilon <1$ and $r>1,$ we have $\eta ^{1-r}\geq 1.$\ \par 
\ \par 
\quad \quad 	Now for $\epsilon >0$ with $\eta _{\epsilon }(z):={\11}_{\Omega
 _{1}}(z)+\epsilon {\11}_{\Omega \backslash \Omega _{1}}(z),$
 let $u_{\epsilon }\in L^{r}(\Omega ,\eta _{\epsilon }^{1-r})$
 be the previous solution, then\ \par 
\vspace{5pt} \ \ \ \ \ \ \ \ \ \ \ \  $\displaystyle {\left\Vert{u_{\epsilon
 }}\right\Vert}_{L^{r}(\Omega ,\eta _{\epsilon }^{1-r})}^{r}\leq
 \int_{\Omega }{\left\vert{u_{\epsilon }}\right\vert ^{r}\eta
 ^{1-r}dm}\leq C^{r}{\left\Vert{\omega }\right\Vert}_{L^{r}(\Omega
 )}^{r}.$\vspace{5pt}  \ \par 
Replacing $\eta $ by its value we get\ \par 
\vspace{5pt} \ \ \ \ \ \ \ \ \ \ \ \  $\displaystyle \int_{\Omega
 _{1}}{\left\vert{u_{\epsilon }}\right\vert ^{r}dm}+\ \int_{\Omega
 \backslash \Omega _{1}}{\left\vert{u_{\epsilon }}\right\vert
 ^{r}\epsilon ^{1-r}dm}\leq C^{r}{\left\Vert{\omega }\right\Vert}_{L^{r}(\Omega
 )}^{r}\Rightarrow $\vspace{5pt}  \ \par 
\vspace{5pt} \ \ \ \ \ \ \ \ \ \ \ \ \ \ \ \ \ \ \ \ \ \ \ \
  $\displaystyle \Rightarrow \int_{\Omega \backslash \Omega _{1}}{\left\vert{u_{\epsilon
 }}\right\vert ^{r}\epsilon ^{1-r}dm}\leq C^{r}{\left\Vert{\omega
 }\right\Vert}_{L^{r}(\Omega )}^{r}$\vspace{5pt}  \ \par 
hence\ \par 
\vspace{5pt} \ \ \ \ \ \ \ \ \ \ \ \  $\displaystyle \int_{\Omega
 \backslash \Omega _{1}}{\left\vert{u_{\epsilon }}\right\vert
 ^{r}dm}\leq C^{r}\epsilon ^{r-1}{\left\Vert{\omega }\right\Vert}_{L^{r}(\Omega
 )}^{r}.$\vspace{5pt}  \ \par 
Because $C$ and the norm of $\omega $ are independent of $\epsilon
 ,$ we have that $\ {\left\Vert{u_{\epsilon }}\right\Vert}_{L^{r}(\Omega
 )}$ is uniformly bounded and $r>1$ implies that $\displaystyle
 L_{p,q-1}^{r}(\Omega )$ is a dual by Lemma~\ref{AG25}, hence
 there is a sub-sequence $\lbrace u_{\epsilon _{k}}\rbrace _{k\in
 {\mathbb{N}}}$ of $\lbrace u_{\epsilon }\rbrace $ which converges
 weakly, when $\epsilon _{k}\rightarrow 0,$ to a $(p,q-1)$ form
 $u$ in $L_{p,q-1}^{r}(\Omega ),$ still with $\ {\left\Vert{u}\right\Vert}_{L_{p,q-1}^{r}(\Omega
 )}\leq C{\left\Vert{\omega }\right\Vert}_{L_{p,q}^{r}(\Omega
 )}.$ Let us write $\displaystyle u_{k}:=u_{\epsilon _{k}}.$\ \par 
\quad To see that this form $u$ is $0\ a.e.$ on $\Omega \backslash
 \Omega _{1}$ let us write the weak convergence:\ \par 
\vspace{5pt} \ \ \ \ \ \ \ \ \ \ \ \  $\displaystyle \forall
 \alpha \in L_{p,q-1}^{r'}(\Omega ),\ {\left\langle{u_{k},\alpha
 }\right\rangle}=\int_{\Omega }{u_{k}\wedge {\overline{*\alpha
 }}}\rightarrow {\left\langle{u,\alpha }\right\rangle}=\int_{\Omega
 }{u\wedge {\overline{*\alpha }}}.$\vspace{5pt}  \ \par 
As usual take $\displaystyle \alpha :=\frac{u}{\left\vert{u}\right\vert
 }{\11}_{E}$ where $E:=\lbrace \left\vert{u}\right\vert >0\rbrace
 \cap (\Omega \backslash \Omega _{1})$ then we get\ \par 
\vspace{5pt} \ \ \ \ \ \ \ \ \ \ \ \  $\displaystyle \int_{\Omega
 }{u\wedge {\overline{*\alpha }}}=\int_{E}{\left\vert{u}\right\vert
 dm}=\lim _{k\rightarrow \infty }\int_{\Omega }{u_{k}\wedge {\overline{*\alpha
 }}}=\lim _{k\rightarrow \infty }\int_{E}{\frac{u_{k}\wedge {\overline{*u}}}{\left\vert{u}\right\vert
 }}.$\vspace{5pt}  \ \par 
Now we have, by H\"older inequalities:\ \par 
\vspace{5pt} \ \ \ \ \ \ \ \ \ \ \ \  $\displaystyle \left\vert{\int_{E}{\frac{u_{k}\wedge
 {\overline{*u}}}{\left\vert{u}\right\vert }}}\right\vert \leq
 {\left\Vert{u_{k}}\right\Vert}_{L^{r}(E)}{\left\Vert{{\11}_{E}}\right\Vert}_{L^{r'}(E)}.$\vspace{5pt}
  \ \par 
But\ \par 
\vspace{5pt} \ \ \ \ \ \ \ \ \ \ \ \  $\displaystyle {\left\Vert{u_{k}}\right\Vert}_{L^{r}(E)}^{r}\leq
 \int_{\Omega \backslash \Omega _{1}}{\left\vert{u_{k}}\right\vert
 ^{r}dm}\leq (\epsilon _{k})^{r-1}C{\left\Vert{\omega }\right\Vert}_{L^{r}(\Omega
 )}\rightarrow 0,\ k\rightarrow \infty $\vspace{5pt}  \ \par 
and ${\left\Vert{{\11}_{E}}\right\Vert}_{L^{r'}(E)}=(m(E))^{1/r'}.$\ \par 
Hence\ \par 
\vspace{5pt} \ \ \ \ \ \ \ \ \ \ \ \  $\displaystyle \ \left\vert{\
 \int_{E}{\left\vert{u}\right\vert dm}}\right\vert =\lim _{k\rightarrow
 \infty }\int_{E}{\frac{u_{k}\wedge {\overline{*u}}}{\left\vert{u}\right\vert
 }}\leq \ $\vspace{5pt}  \ \par 
\vspace{5pt} \ \ \ \ \ \ \ \ \ \ \ \ \ \ \ \ \ \ \ \ \ \ \ \
  $\displaystyle \leq \lim _{k\rightarrow \infty }C^{r}(m(E))^{1/r'}(\epsilon
 _{k})^{r-1}{\left\Vert{\omega }\right\Vert}_{L^{r}(\Omega )}^{r}=0,$\vspace{5pt}
  \ \par 
so $\int_{E}{\left\vert{u}\right\vert dm}=0$ which implies $m(E)=0$
 because on $E,\ \left\vert{u}\right\vert >0.$\ \par 
\ \par 
\quad Hence we get that the form $u$ is $0\ a.e.$ on $\Omega \backslash
 \Omega _{1}.$\ \par 
\quad \quad 	So we proved\ \par 
\vspace{5pt} \ \ \ \ \ \ \ \ \ \ \ \  $\displaystyle \forall
 \varphi \in {\mathcal{D}}_{n-p,n-q}(\Omega ),\ (-1)^{p+q-1}\ll
 \varphi ,\omega \gg =\ll \bar \partial \varphi ,u_{\epsilon
 }\gg \rightarrow \ll \bar \partial \varphi ,u\gg $\vspace{5pt}  \ \par 
\vspace{5pt} \ \ \ \ \ \ \ \ \ \ \ \ \ \ \ \ \ \ \ \ \ \ \ \
  $\displaystyle \Rightarrow \ll \bar \partial \varphi ,u\gg
 =(-1)^{p+q-1}\ll \varphi ,\omega \gg $\vspace{5pt}  \ \par 
hence $\bar \partial u=\omega $ in the sense of distributions.\ \par 
\quad The proof is complete. $\hfill\blacksquare $\ \par 

\begin{rem}
As in remark~\ref{AndreGrau2815} if $\Omega $ is Stein for $q=n$
 instead of asking $\omega \perp {\mathcal{H}}_{p}(\Omega _{2})$
 we need just $\omega \perp {\mathcal{H}}_{p}(\Omega ).$
\end{rem}

\begin{rem}
The condition of orthogonality to $\displaystyle {\mathcal{H}}_{p}(V)$
 in the case $q=n$ is necessary: suppose there is a $(p,n-1)$
 current $u$ such that $\bar \partial u=\omega $ and $u$ with
 compact support in an open set $V\subset \Omega ,$ then if $\displaystyle
 h\in {\mathcal{H}}_{p}(V),$ we have\par 
\vspace{5pt} \ \ \ \ \ \ \ \ \ \ \ \  $\displaystyle h\in {\mathcal{H}}_{p}(V),\
 \ll \omega ,h\gg =\ll \bar \partial u,h\gg =(-1)^{n+p}\ll u,\bar
 \partial h\gg =0,$\vspace{5pt}  \par 
because, $u$ being compactly supported, there is no boundary term and\par 
\vspace{5pt} \ \ \ \ \ \ \ \ \ \ \ \  $\displaystyle \ll \bar
 \partial u,h\gg =(-1)^{n+p}\ll u,\bar \partial h\gg .$\vspace{5pt}  \par 
This kind of condition was already seen for extension of CR 
 functions, see~\cite{AmExtCR91} and the references therein.
\end{rem}

\subsection{Finer control of the support.~\label{AndreGrau1913}}
\quad \quad 	Here we shall get a better control on the support of a solution.\ \par 

\begin{thm}
~\label{AndreGrau1911}Let $\Omega $ be a weakly $r'$ regular
 domain in a Stein manifold $X.$\par 
\quad Suppose the $(p,q)$ form $\omega $ is in $L^{r,c}(\Omega ,dm),\
 \bar \partial \omega =0,$ if $q<n,$ and $\omega \perp {\mathcal{H}}_{p}(V)$
 for any $V$ such that $\displaystyle \Supp \omega \subset V,$
 if $q=n,$ with  $\Supp \omega \subset \Omega \backslash C,$
 where $C$ is a weakly $r$ regular domain.\par 
\quad For any open relatively compact set $U$ in $C,$ there is a $u\in
 L^{r,c}(\Omega ,dm)$ such that $\bar \partial u=\omega $ and
 with support in $\Omega \backslash \bar U,$ provided that $q\geq 2.$
\end{thm}
\quad \quad 	Proof.\ \par 
Let $\omega $ be a $(p,q)$ form with compact support in $\Omega
 \backslash C$ then there is a  $v\in L^{r}_{p,q-1}(\Omega ),\
 \bar \partial v=\omega ,$ with compact support in $\Omega ,$
 by theorem~\ref{AndGrau84} or, if $\Omega $ is a polydisc in
 ${\mathbb{C}}^{n}$ and if $\omega \in {\mathcal{W}}^{r}_{q}(\Omega
 ),$ by the theorem in~\cite{AmMong12}.\ \par 
\quad Because $\omega $ has compact support outside $C$ we have $\omega
 =0$ in $C;$ this means that $\bar \partial v=0$ in $C.$ Because
 $C$ is weakly $r$ regular and $q\geq 2,$ we have\ \par 
\vspace{5pt} \ \ \ \ \ \ \ \ \ \ \ \  $\displaystyle \exists
 C'\subset C,\ C'\supset \bar U,\ \exists h\in L^{r}_{p,q-2}(C')\
 s.t.\ \bar \partial h=v$\vspace{5pt}   in $C'.$\vspace{5pt}  \ \par 
Let $\chi $ be a smooth function such that $\chi =1$ in $U$ and
 $\chi =0$ near $\partial C';$ then set	 	$\displaystyle u:=v-\bar
 \partial (\chi h).$\ \par 
We have that $u=v-\chi \bar \partial h-\bar \partial \chi \wedge
 h=v-\chi v-\bar \partial \chi \wedge h$ hence $u$ is in $L^{r}(\Omega
 );$ moreover $u=0$ in $\bar U$ because $\chi =1$ in $U$ hence
 $\bar \partial \chi =0$ there. Finally $\bar \partial u=\bar
 \partial v-\bar \partial ^{2}(\chi h)=\omega $ and we are done.
 $\hfill\blacksquare $\ \par 
\ \par 
\quad \quad 	If $\Omega $ and $C$ are, for instance, pseudo-convex in ${\mathbb{C}}^{n}$
 then $\Omega \backslash C$ is no longer pseudo-convex in general,
 so this theorem improves actually the control of the support.\ \par 

\begin{rem}
The correcting function $h$ is given by kernels in the case of
 Stein domains, hence it is linear; if the primitive solution
 $v$ is also linear in $\omega ,$ then the solution $u$ is linear
 too. This is the case in ${\mathbb{C}}^{n}$ with the solution
 given in~\cite{AmMong12}.
\end{rem}
\quad \quad 	This theorem cannot be true for $q=1$ as shown by the following example:\ \par 
take a holomorphic function $\varphi $ in the open unit ball
 $B(0,1)$ in ${\mathbb{C}}^{n}$ such that it extends to no open
 ball of center $0$ and radius $>1.$ For instance $\varphi (z):=\exp
 {\left({-\frac{z_{1}+1}{z_{1}-1}}\right)}.$ Take  $R<1,$ then
 $\varphi $ is ${\mathcal{C}}^{\infty }(\bar B(0,R))$ hence by
 a theorem of Whitney $\varphi $ extends ${\mathcal{C}}^{\infty
 }$ to ${\mathbb{C}}^{n};$ call  $\varphi _{R}$ this extension.
 Let $\chi \in {\mathcal{C}}^{\infty }_{c}(B(0,2))$ such that
 $\chi =1$ in the ball $B(0,3/2)$ and consider the $(0,1)$ form
 $\omega :=\bar \partial (\chi \varphi _{R}).$ Then $\Supp \omega
 \subset B(0,2)\backslash B(0,R),\ \omega $ is $\bar \partial
 $ closed and is ${\mathcal{C}}^{\infty }$ hence in $L^{r}_{0,1}(B(0,2)).$
 Moreover $B(0,R)$ is strictly pseudo-convex hence $r'$ regular,
 but there is no function $u$ such that $\bar \partial u=\omega
 $ and $u$ zero near the origin because any solution $u$ will
 be C.R. on $\partial B(0,R)$ and by Hartog's phenomenon will
 extends holomorphically to $B(0,R),$ hence cannot be identically
  null near $0.$\ \par 
\ \par 
\quad \quad 	Never the less in the case $q=1,$ we have:\ \par 

\begin{thm}
~\label{AndreGrau1912} Let $\Omega $ be a weakly $r'$ regular
 domain in a Stein manifold $X.$\par 
Then for any $(p,1)$ form $\omega $ in $L^{r,c}(\Omega ),\ \bar
 \partial \omega =0,$ with support in $\Omega _{1}\backslash
 C$ where $\Omega _{1}$ is a weak $r'$ regular domain in  $\Omega
 $ and $C$ is a domain such that $C\subset \Omega $ and $C\backslash
 \Omega _{1}\neq \emptyset ;$ there is a $\displaystyle u\in
 L^{r,c}(\Omega )$ such that $\bar \partial u=\omega $ and with
 support in $\Omega \backslash C.$
\end{thm}
\quad \quad 	Proof.\ \par 
There is $u\in L^{r}_{p,0}(\Omega _{1})$ such that $\bar \partial
 u=\omega $ with compact support in $\Omega _{1},$ by theorem~\ref{AndGrau84}.
 Then $\bar \partial u=0$ in $C$ hence $u$ is locally holomorphic
 in $C.$ Because $C\backslash \Omega _{1}\neq \emptyset ,$ there
 is an open set in $C\backslash \Omega _{1}\subset \Omega \backslash
 \Omega _{1}$ where $u$ is $0$ and holomorphic, hence $u$ is
 identically $0$ in $C,\ C$ being connected. $\hfill\blacksquare $\ \par 

\begin{rem}
If there is a $u\in L^{r,c}_{p,0}(\Omega _{1})$ which is $0$
 in $C,$ we have\par 
\vspace{5pt} \ \ \ \ \ \ \ \ \ \ \ \  $\displaystyle \forall
 h\in L^{r'}_{n-p,n-1}(C)::\Supp \bar \partial h\subset C,\ 0=\ll
 u,\bar \partial h\gg =\ll \omega ,h\gg ,$\vspace{5pt}  \par 
hence the necessary condition:\par 
\vspace{5pt} \ \ \ \ \ \ \ \ \ \ \ \  $\displaystyle \forall
 h\in L^{r'}_{n-p,n-1}(C)::\Supp \bar \partial h\subset C,\ \ll
 \omega ,h\gg =0.$\vspace{5pt}  
\end{rem}
\quad We proved in~\cite{AmMong12}:\ \par 

\begin{thm}
Let $f\in {\mathcal{O}}(\bar {\mathbb{D}}^{n})$ be a holomorphic
 function in a neighborhood of the closed unit polydisc in ${\mathbb{C}}^{n}$
 and set $Z:=f^{-1}(0).$ Then for any $(0,q)$ form $\omega $
 in $L^{r}({\mathbb{D}}^{n}\backslash Z)\cap {\mathcal{W}}^{r}_{q}(\Omega
 ),\ \bar \partial \omega =0,$ with compact support in ${\mathbb{D}}^{n}\backslash
 Z,$ for any $k\in {\mathbb{N}},$ we can find a $(0,q-1)$-form
 $\beta \in L^{r,c}({\mathbb{D}}^{n})$ such that $\bar \partial
 (f^{k}\beta )=\omega .$ Equivalently we can find a $(0,q-1)$-form
  $\eta =f^{k}\beta $ such that $\eta \in L^{r,c}({\mathbb{D}}^{n}),\
 \eta $ is $0$ on $Z$ up to order $k$ and $\bar \partial \eta =\omega .$
\end{thm}
\quad And by Remark 6.3 of this paper, the solutions are given by a
 \emph{bounded linear operator.}\ \par 
\ \par 
\quad The following corollary will generalise strongly this result
 but at the price that \emph{we have not the linearity, nor even
 the constructivity of the solution.}\ \par 

\begin{cor}
Let  $\Omega $ be a Stein manifold.  Let $f$ be a holomorphic
 function in $\Omega $ and set $Z:=f^{-1}(0).$  Then for any
 $(p,q)$ form $\omega $ in $L^{r,c}(\Omega \backslash Z),\ \bar
 \partial \omega =0,$ if $1\leq q<n,$ and $\omega \perp {\mathcal{H}}_{p}(\
 \Omega \backslash Z)$ if $q=n,$ there is a $(p,q-1)$ form $\displaystyle
 u\in L^{r}(\Omega \backslash Z)$ such that $\bar \partial u=\omega
 $ and $u$ has its support still in $\Omega \backslash Z.$
\end{cor}
\quad Proof.\ \par 
We first show that $\displaystyle \Omega \backslash Z$ is Stein.
 Because $f\neq 0$ in $\displaystyle \Omega \backslash Z$ we
 have that $\varphi :=\frac{1}{\left\vert{f}\right\vert ^{2}}$
 is plurisubharmonic in $\displaystyle \Omega \backslash Z$ and
 ${\mathcal{C}}^{\infty }(\Omega \backslash Z).$  Because $\Omega
 $ is Stein we have, by Theorem 5.1.6 of H\"ormander~\cite{Hormander73},
 a strictly plurisubharmonic exhausting function $\rho $ in ${\mathcal{C}}^{\infty
 }(\Omega ).$ Now the function $\gamma :=\varphi +\rho $ is still
 strictly plurisubharmonic and ${\mathcal{C}}^{\infty }$ in $\displaystyle
 \Omega \backslash Z.$ Now we shall prove:\ \par 
\vspace{5pt} \ \ \ \ \ \ \ \ \ \ \ \  $\displaystyle \forall
 \alpha \in {\mathbb{R}},\ K_{\alpha }:=\lbrace z\in \Omega \backslash
 Z::\gamma (z)<\alpha \rbrace $\vspace{5pt}   is relatively compact
 in $\Omega \backslash Z.$\vspace{5pt}  \ \par 
We have $\rho (z)<\alpha -\varphi (z)<\alpha $ on $K_{\alpha
 }$ because $\varphi (z)\geq 0,$ hence, because $\rho $ is exhaustive
 in $\Omega ,$ we have that $K_{\alpha }$ is contained in a compact
 set $F$ in $\Omega .$ So on $F,$ hence on $K_{\alpha },$ we
 have that $\rho (z)\geq A>-\infty $ because $\rho $ is continuous.\ \par 
We also have $\varphi (z)<\alpha -\rho (z)$ on $K_{\alpha }$
 i.e. $\left\vert{f(z)}\right\vert ^{2}>\frac{1}{\alpha -\rho
 (z)}.$ So, on the set $K_{\alpha },\ \alpha >\rho (z)\geq A>-\infty
 ,$ hence $\displaystyle \left\vert{f(z)}\right\vert >\frac{1}{\alpha
 -A}>0$ on $K_{\alpha },$ so $K_{\alpha }$ is far away from $Z,$
 hence $K_{\alpha }$ is relatively compact in $\displaystyle
 \Omega \backslash Z.$\ \par 
\quad So we can apply ~\cite[Theorem 5.2.10, p. 127]{Hormander73} to
 get that $\displaystyle \Omega \backslash Z$ is a Stein manifold.\ \par 
\quad Now we are in position to apply Theorem~\ref{AndGrau84}. Let
 $\omega $ be a $(p,q)$ form in $L^{r,c}(\Omega \backslash Z),\
 \bar \partial \omega =0,$ if $1\leq q<n,$ and $\omega \perp
 {\mathcal{H}}_{p}(\ \Omega \backslash Z)$ if $q=n,$ Theorem~\ref{AndGrau84}
 gives a $(p,q-1)$ form $\displaystyle u\in L^{r}(\Omega \backslash
 Z)$ such that $\bar \partial u=\omega $ and $u$ has its compact
 support in $\Omega \backslash Z.$\ \par 
The proof is complete. $\hfill\blacksquare $\ \par 

\begin{rem}
This leaves open the question to have a linear (or a constructive)
 solution to this problem even in the case of the polydisc.
\end{rem}

\section{Appendix}
\quad Here we shall prove certainly known results on the duality $L^{r}-L^{r'}$
 for $(p,q)$-forms in a complex manifold $X.$ Because I was unable
 to find precise references for them, I prove them here.\ \par 
\ \par 
Recall we have a pointwise scalar product and a pointwise modulus
 for $(p,q)$-forms in $X$:\ \par 
\vspace{5pt} \ \ \ \ \ \ \ \ \ \ \ \  $\displaystyle (\alpha
 ,\beta )dm:=\alpha \wedge {\overline{*\beta }};\ \ \left\vert{\alpha
 }\right\vert ^{2}dm:=\alpha \wedge {\overline{*\alpha .}}$\vspace{5pt}  \ \par 
By the Cauchy-Schwarz inequality for scalar products we get:\ \par 
\vspace{5pt} \ \ \ \ \ \ \ \ \ \ \ \  $\displaystyle \forall
 x\in X,\ \left\vert{(\alpha ,\beta )(x)}\right\vert \leq \left\vert{\alpha
 (x)}\right\vert \left\vert{\beta (x)}\right\vert .$\vspace{5pt}  \ \par 
\quad This gives H\"older inequalities for $(p,q)$-forms:\ \par 

\begin{lem}
~\label{AG24}(H\"older inequalities) Let $\alpha \in L^{r}_{p,q}(\Omega
 )$ and $\displaystyle \beta \in L^{r'}_{p,q}(\Omega ).$ We have\par 
\vspace{5pt} \ \ \ \ \ \ \ \ \ \ \ \  $\displaystyle \left\vert{{\left\langle{\alpha
 ,\beta }\right\rangle}}\right\vert \leq {\left\Vert{\alpha }\right\Vert}_{L^{r}(\Omega
 )}{\left\Vert{\beta }\right\Vert}_{L^{r'}(\Omega )}.$\vspace{5pt}  
\end{lem}
\quad Proof.\ \par 
We start with ${\left\langle{\alpha ,\beta }\right\rangle}=\int_{\Omega
 }{(\alpha ,\beta )(x)dm(x)}$ hence\ \par 
\vspace{5pt} \ \ \ \ \ \ \ \ \ \ \ \  $\displaystyle \left\vert{{\left\langle{\alpha
 ,\beta }\right\rangle}}\right\vert \leq \int_{\Omega }{\left\vert{(\alpha
 ,\beta )(x)}\right\vert dm}\leq \int_{\Omega }{\left\vert{\alpha
 (x)}\right\vert \left\vert{\beta (x)}\right\vert dm(x)}.$\vspace{5pt}  \ \par 
By the usual H\"older inequalities for functions we get\ \par 
\vspace{5pt} \ \ \ \ \ \ \ \ \ \ \ \  $\displaystyle \int_{\Omega
 }{\left\vert{\alpha (x)}\right\vert \left\vert{\beta (x)}\right\vert
 dm(x)}\leq {\left({\int_{\Omega }{\left\vert{\alpha (x)}\right\vert
 ^{r}dm}}\right)}^{1/r}{\left({\int_{\Omega }{\left\vert{\beta
 (x)}\right\vert ^{r'}dm}}\right)}^{1/r'}$\vspace{5pt}  \ \par 
which ends the proof of the lemma. $\hfill\blacksquare $\ \par 

\begin{lem}
~\label{AG21}Let $\alpha \in L_{p,q}^{r}(\Omega )$ then\par 
\vspace{5pt} \ \ \ \ \ \ \ \ \ \ \ \  $\displaystyle {\left\Vert{\alpha
 }\right\Vert}_{L_{p,q}^{r}(\Omega )}=\sup _{\beta \in L_{p,q}^{r'}(\Omega
 ),\ \beta \neq 0}\frac{\left\vert{{\left\langle{\alpha ,\beta
 }\right\rangle}}\right\vert }{{\left\Vert{\beta }\right\Vert}_{L^{r'}(\Omega
 )}}.$\vspace{5pt}  
\end{lem}
\quad Proof.\ \par 
We choose $\beta :=\alpha \left\vert{\alpha }\right\vert ^{r-2},$ then:\ \par 
\vspace{5pt} \ \ \ \ \ \ \ \ \ \ \ \  $\displaystyle \left\vert{\beta
 }\right\vert ^{r'}=\left\vert{\alpha }\right\vert ^{r'(r-1)}=\left\vert{\alpha
 }\right\vert ^{r}\Rightarrow {\left\Vert{\beta }\right\Vert}_{L^{r'}(\Omega
 )}^{r'}={\left\Vert{\alpha }\right\Vert}_{L^{r}(\Omega )}^{r}.$\vspace{5pt}
  \ \par 
Hence\ \par 
\vspace{5pt} \ \ \ \ \ \ \ \ \ \ \ \  $\displaystyle {\left\langle{\alpha
 ,\beta }\right\rangle}={\left\langle{\alpha ,\alpha \left\vert{\alpha
 }\right\vert ^{r-2}}\right\rangle}=\int_{\Omega }{(\alpha ,\alpha
 )\left\vert{\alpha }\right\vert ^{r-2}dm}={\left\Vert{\alpha
 }\right\Vert}_{L^{r}(\Omega )}^{r}.$\vspace{5pt}  \ \par 
On the other hand we have\ \par 
\vspace{5pt} \ \ \ \ \ \ \ \ \ \ \ \  $\displaystyle {\left\Vert{\beta
 }\right\Vert}_{L^{r'}(\Omega )}={\left\Vert{\alpha }\right\Vert}_{L^{r}(\Omega
 )}^{r/r'}={\left\Vert{\alpha }\right\Vert}_{L^{r}(\Omega )}^{r-1},$\vspace{5pt}
  \ \par 
so\ \par 
\vspace{5pt} \ \ \ \ \ \ \ \ \ \ \ \  $\displaystyle {\left\Vert{\alpha
 }\right\Vert}_{L^{r}(\Omega )}{\times}{\left\Vert{\beta }\right\Vert}_{L^{r'}(\Omega
 )}={\left\Vert{\alpha }\right\Vert}_{L^{r}(\Omega )}^{r}={\left\langle{\alpha
 ,\beta }\right\rangle}.$\vspace{5pt}  \ \par 
Hence $\displaystyle {\left\Vert{\alpha }\right\Vert}_{L^{r}(\Omega
 )}=\frac{\left\vert{{\left\langle{\alpha ,\beta }\right\rangle}}\right\vert
 }{{\left\Vert{\beta }\right\Vert}_{L^{r'}(\Omega )}}.$\ \par 
A fortiori for any choice of $\beta $:\ \par 
\vspace{5pt} \ \ \ \ \ \ \ \ \ \ \ \  $\displaystyle {\left\Vert{\alpha
 }\right\Vert}_{L^{r}(\Omega )}\leq \sup _{\beta \in L^{r'}(\Omega
 )}\frac{\left\vert{{\left\langle{\alpha ,\beta }\right\rangle}}\right\vert
 }{{\left\Vert{\beta }\right\Vert}_{L^{r'}(\Omega )}}.$\vspace{5pt}  \ \par 
To prove the other direction, we use the H\"older inequalities,
 Lemma~\ref{AG24}:\ \par 
\vspace{5pt} \ \ \ \ \ \ \ \ \ \ \ \  $\displaystyle \forall
 \beta \in L_{p,q}^{r'}(\Omega ),\ \frac{\left\vert{{\left\langle{\alpha
 ,\beta }\right\rangle}}\right\vert }{{\left\Vert{\beta }\right\Vert}_{L^{r'}(\Omega
 )}}\leq {\left\Vert{\alpha }\right\Vert}_{L^{r}(\Omega )}.$\vspace{5pt}
  \ \par 
The proof is complete. $\hfill\blacksquare $\ \par 
\ \par 
\quad Now we are in position to state:\ \par 

\begin{lem}
~\label{AG25}The dual space of the Banach space $\displaystyle
 L_{p,q}^{r}(\Omega )$ is $\displaystyle L_{n-p,n-q}^{r'}(\Omega ).$
\end{lem}
\quad Proof.\ \par 
Suppose first that $u\in L_{n-p,n-q}^{r'}(\Omega ).$ Then consider:\ \par 
\vspace{5pt} \ \ \ \ \ \ \ \ \ \ \ \  $\displaystyle \forall
 \alpha \in L_{p,q}^{r}(\Omega ),\ {\mathcal{L}}(\alpha ):=\int_{\Omega
 }{\alpha \wedge u}={\left\langle{\alpha ,{\overline{*u}}}\right\rangle}.$\vspace{5pt}
  \ \par 
This is a linear form on $\displaystyle L_{p,q}^{r}(\Omega )$
 and its norm, by definition, is\ \par 
\vspace{5pt} \ \ \ \ \ \ \ \ \ \ \ \  $\displaystyle {\left\Vert{{\mathcal{L}}}\right\Vert}=\sup
 _{\alpha \in L^{r}(\Omega )}\frac{\left\vert{{\left\langle{\alpha
 ,{\overline{*u}}}\right\rangle}}\right\vert }{{\left\Vert{\alpha
 }\right\Vert}_{L^{r}(\Omega )}}.$\vspace{5pt}  \ \par 
By use of Lemma~\ref{AG21} we get\ \par 
\vspace{5pt} \ \ \ \ \ \ \ \ \ \ \ \  $\displaystyle {\left\Vert{{\mathcal{L}}}\right\Vert}={\left\Vert{{\overline{*u}}}\right\Vert}_{L^{r'}_{p,q}(\Omega
 )}={\left\Vert{u}\right\Vert}_{L^{r'}_{n-p,n-q}(\Omega )}.$\vspace{5pt}
  \ \par 
So we have ${\left({L_{p,q}^{r}(\Omega )}\right)}^{*}\supset
 L_{n-p,n-q}^{r'}(\Omega )$ with the same norm.\ \par 
\ \par 
\quad Conversely take a continuous linear form ${\mathcal{L}}$ on $\displaystyle
 L_{p,q}^{r}(\Omega ).$ We have, again by definition, that:\ \par 
\vspace{5pt} \ \ \ \ \ \ \ \ \ \ \ \  $\displaystyle {\left\Vert{{\mathcal{L}}}\right\Vert}=\sup
 _{\alpha \in L^{r}(\Omega )}\frac{\left\vert{{\mathcal{L}}(\alpha
 )}\right\vert }{{\left\Vert{\alpha }\right\Vert}_{L^{r}(\Omega
 )}}.$\vspace{5pt}  \ \par 
Because ${\mathcal{D}}_{p,q}(\Omega )\subset L_{p,q}^{r}(\Omega
 ),$  ${\mathcal{L}}$ is a continuous linear form on $\displaystyle
 {\mathcal{D}}_{p,q}(\Omega ),$ hence, by definition, ${\mathcal{L}}$
 can be represented by a $(n-p,n-q)$-current $u.$ So we have:\ \par 
\vspace{5pt} \ \ \ \ \ \ \ \ \ \ \ \  $\displaystyle \forall
 \alpha \in {\mathcal{D}}_{p,q}(\Omega ),\ {\mathcal{L}}(\alpha
 ):=\int_{\Omega }{\alpha \wedge u}={\left\langle{\alpha ,{\overline{*u}}}\right\rangle}.$\vspace{5pt}
  \ \par 
Moreover we have, by Lemma~\ref{AG21},\ \par 
\vspace{5pt} \ \ \ \ \ \ \ \ \ \ \ \  $\displaystyle {\left\Vert{{\mathcal{L}}}\right\Vert}=\sup
 _{\alpha \in {\mathcal{D}}_{p,q}(\Omega )}\frac{\left\vert{{\left\langle{\alpha
 ,*\bar u}\right\rangle}}\right\vert }{{\left\Vert{\alpha }\right\Vert}_{L^{r}(\Omega
 )}}={\left\Vert{*u}\right\Vert}_{L^{r'}(\Omega )}$\vspace{5pt}  \ \par 
because $\displaystyle {\mathcal{D}}_{p,q}(\Omega )$ is dense
 in $\displaystyle L_{p,q}^{r}(\Omega ).$ So we proved $\displaystyle
 {\left({L_{p,q}^{r}(\Omega )}\right)}^{*}\subset L_{n-p,n-q}^{r'}(\Omega
 )$ with the same norm.\ \par 
The proof is complete. $\hfill\blacksquare $\ \par 
\ \par 

\bibliographystyle{/usr/local/texlive/2017/texmf-dist/bibtex/bst/base/apalike}

\end{document}